\documentclass[a4paper,12pt,eqno]{article}
\usepackage{theorem}
\usepackage{latexsym,amssymb,amsfonts,amsmath}
\usepackage{graphicx}
\newtheorem{thm}{Theorem}[section]

\newtheorem{cor}[thm]{Corollary}
\newtheorem{pro}[thm]{Proposition}

\newtheorem{prob}[thm]{Problem} 

\theoremheaderfont{\scshape}

\title{{\Large {\bf A zeta function related to the transition matrix of 
the discrete-time quantum walk on a graph}}
\author{
{\small Norio KONNO}\\
{\scriptsize Department of Applied Mathematics, 
Faculty of Engineering, Yokohama National University}\\
{\scriptsize Hodogaya, Yokohama 240-8501, Japan}\\
{\small Iwao SATO} \\ 
{\scriptsize Oyama National College of Technology} \\ 
{\scriptsize Oyama, Tochigi 323-0806, Japan} \\
{\small Etsuo SEGAWA} \\
{\scriptsize Graduate School of Education Center, Yokohama National University} \\
{\scriptsize Hodogaya, Yokohama, 240-8501, Japan}.}
}
\date{\empty }
\begin{document}
\maketitle

\par\noindent
\begin{small}
\par\noindent
{\bf Abstract}.  We present the structure theorem for the positive support of the cube of 
the Grover transition matrix of the discrete-time quantum walk (the Grover walk) on a general graph $G$ 
under same condition. 
Thus, we introduce a zeta function on the positive support of the cube of 
the Grover transition matrix of $G$, and present its Euler product and its determinant expression. 
As a corollary, we give the characteristic polynomial for the positive support of the cube of 
the Grover transition matrix of a regular graph, and so obtain its spectra. 
Finally, we present the poles and  the radius of the convergence of this zeta function.   
\footnote[0]{
{\it Abbr. title:} The transition matrix of a quantum walk on a graph 
}
\end{small}

\setcounter{equation}{0}
\section{Introduction}
As a quantum counterpart of the classical random walk, the quantum walk (QW) has recently attracted much attention for various fields. 
There are two types of QWs. One is the discrete-time walk and the other is the continuous-time one. 
The discrete-time QW in one dimension (1D) was intensively studied by Ambainis et al. \cite{AmbainisEtAl2001}. 
One of the most striking properties of the 1D QW is the spreading property of the walker. 
The standard deviation of the position grows linearly in time, quadratically faster than classical random walk. 
The review and book on QWs are Kempe \cite{Kempe2003}, Konno \cite{Konno2008b}, for examples.

Recently, quantum walks of graphs were applied in graph isomorphism problems(see \cite{ShiauETAL},\cite{EmmsETAL2009},
\cite{GambelETAL},\cite{EmmsETAL2006}). 
Emms et al. \cite{EmmsETAL2006} treated spectra of the transition matrix(Grover transition matrix), its positive support and the 
positive support of its square of the discrete-time Grover walk on a graph, and showed that the third power of the transition matrix 
outperforms the graph spectra methods in distinguishing strongly regular graphs. 
Godsil and Guo \cite{GG2010} gave new proofs of the results of Emms et al. \cite{EmmsETAL2006}.  
Higuchi, Konno, Sato and Segwa \cite{HKSS} introduced a generalized Szegedy transition matrix as a generalization 
of the transition matrix of the Szegedy walk on a graph, and obtained the results of Emms et al. \cite{EmmsETAL2006} from 
the characteristic polynomial of the generalized Szegedy transition matrix. 
Furthermore, they \cite{HKSS} decided the structure of the positive support of the cube of 
the Grover transition matrix of a regular graph under some condition.  

It is turned out that the zeta function of a graph is closely related to the discrete-time quantum walk on a graph. 
Already, the Ihara zeta function of a graph obtained various success related to graph spectra(see \cite{Ihara1966},
\cite{Sunada1986},\cite{Sunada1988},\cite{Hashimoto1989},\cite{Bass1992},\cite{KS2000}). 
Ihara zeta functions of graphs started from Ihara zeta functions of regular 
graphs by Ihara \cite{Ihara1966}. 
Hashimoto \cite{Hashimoto1989} and Bass \cite{Bass1992} generalized Ihara's result on the zeta function of 
a regular graph to an irregular graph, and showed that its reciprocal is again a polynomial. 
Kotani and Sunada \cite{KS2000} presented the distribution of poles and the radius of the convergence of 
the Ihara zeta function.  

Ren et al. \cite{RenETAL} found a interesting relationship between the Ihara zeta function and 
the discrete-time quantum walks on a graph, and showed that the positive support of 
the transition matrix of the discrete-time quantum walk is equal to 
the Perron-Frobenius operator(the edge matrix) related to the Ihara zeta function. 
Konno and Sato \cite{KS2012} presented the characteristic polynomials of the Grover transition matrix and its positive support 
of a graph by using the determinant expressions of the second weighted zeta function and the Ihara zeta function of a graph, 
and so decided spectra of the Grover transition matrix and its positive support from them. 
Higuch, Konno, Sato and Segawa  \cite{HKSS2014} introduced we introduce a zeta function of a graph with respect to the positive support of 
the square of the Grover transition matrix, and present its Euler product, its determinant expression, its poles and its radius of 
the convergence. 
 
The rest of the paper is organized as follows. 
Section 2 give a short survey about the positive support of Grover transition matrix etc, and their spectra. 
Section 3 present a short survey on the Ihara zeta function of a graph. 
In Sect. 4, we present the structure theorem for the positive support of the cube of the Grover transition matrix of 
a general graph $G$ under same condition to that of Higuch, Konno, Sato and Segawa  \cite{HKSS2014}. 
In Sect. 5, we introduce a zeta function on the positive support of the cube of 
the Grover transition matrix of $G$, and present its Euler product and its determinant expression. 
As a corollary, we give the characteristic polynomial for the positive support of the cube of 
the Grover transition matrix of a regular graph, and so obtain its spectra. 
In Sect. 6, we present the poles and  the radius of the convergence of the above zeta function. 

\section{The transition matrix of a quantum walk on a graph} 

Graphs treated here are finite.
Let $G=(V(G),E(G))$ be a connected graph (possibly multiple edges and loops) 
with the set $V(G)$ of vertices and the set $E(G)$ of unoriented edges $uv$ 
joining two vertices $u$ and $v$. 
For $uv \in E(G)$, an arc $(u,v)$ is the oriented edge from $u$ to $v$. 
Set $D(G)= \{ (u,v),(v,u) \mid uv \in E(G) \} $. 
For $e=(u,v) \in D(G)$, set $u=o(e)$ and $v=t(e)$. 
Furthermore, let $e^{-1}=(v,u)$ be the {\em inverse} of $e=(u,v)$. 
The {\em degree} $\deg v = \deg {}_G \  v$ of a vertex $v$ of $G$ is the number of edges incident to $v$. 
For a natural number $k$, a graph $G$ is called {\em $k$-regular } if $\deg {}_G \  v=k$ for each vertex $v$ of $G$. 

A discrete-time quantum walk is a quantum process on a graph whose state vector is governed by 
a matrix called the transition matrix.  
Let $G$ be a connected graph with $n$ vertices and $m$ edges, 
$V(G)= \{ v_1 , \ldots , v_n \} $ and $D(G)= \{ e_1 , \ldots , e_m , 
e^{-1}_1 , \ldots , e^{-1}_m \} $. 
Set $d_j = d_{u_j} = \deg v_j $ for $i=1, \ldots ,n$. 
The the {\em transition matrix} ${\bf U} ={\bf U} (G)=( U_{ef} )_{e,f \in D(G)} $ 
of $G$ is defined by 
\[
U_{ef} =\left\{
\begin{array}{ll}
2/d_{t(f)} (=2/d_{o(e)} ) & \mbox{if $t(f)=o(e)$ and $f \neq e^{-1} $, } \\
2/d_{t(f)} -1 & \mbox{if $f= e^{-1} $, } \\
0 & \mbox{otherwise}
\end{array}
\right. 
\]

We introduce the {\em positive support} ${\bf F}^+ =( F^+_{ij} )$ of 
a real square matrix ${\bf F} =( F_{ij} )$ as follows: 
\[
F^+_{ij} =\left\{
\begin{array}{ll}
1 & \mbox{if $F_{ij} >0$, } \\
0 & \mbox{otherwise}
\end{array}
\right.
\]

Emms et al \cite{EmmsETAL2006} expressed the spectra of the positive support ${\bf U}^+ $ of 
the Grover transition matrix of a regular graph $G$ by means of those of the adjacency matrix ${\bf A} (G)$ of $G$.

\begin{thm}[Emms, Hancock, Severini and Wilson \cite{EmmsETAL2006}]
\label{thm5.1}
Let $G$ be a connected $k$-regular graph with $n$ vertices and $m$ edges, and 
$\delta (G) \geq 2$. 
The positive support ${\bf U}^+ $ has $2n$ eigenvalues of the form 
\[
\lambda = \frac{\lambda {}_A }{2} \pm i \sqrt{k-1- \lambda {}^2_A /4} , 
\]
where $\lambda {}_A $ is an eigenvalue of the matrix ${\bf A} (G)$. 
The remaining $2(m-n)$ eigenvalues of ${\bf U}^+$ are $\pm 1$ with equal multiplicities. 
\end{thm}

Furthermore, Emms et al \cite{EmmsETAL2006} expressed the spectra of the positive support ${\bf U}^+ $ of the square   
of the transition matrix of a regular graph $G$ by means of those of 
the adjacency matrix ${\bf A} (G)$ of $G$.

\begin{thm}[Emms, Hancock, Severini and Wilson \cite{EmmsETAL2006}]
Let $G$ be a connected $k$-regular graph with $n$ vertices and $m$ edges. 
Suppose that $k>2$. 
The positive support $({\bf U}^2 ) ^+ $ has $2n$ eigenvalues of the form 
\[
\lambda = \frac{\lambda {}^2_A -2k+4}{2} \pm i \frac{\lambda {}_A \sqrt{4k-4- \lambda {}^2_A }}{2} , 
\]
where $\lambda {}_A $ is an eigenvalue of the matrix ${\bf A} (G)$. 
The remaining $2(m-n)$ eigenvalues of ${\bf U}^+$ are $2$. 
\end{thm}

Next, we state the structure of $({\bf U}^+ $, $( {\bf U}^2)^+ $ and $( {\bf U}^3)^+ $. 

If the degree of each vertex of $G$ is not less than 2, i.e., $\delta (G) \geq 2$, 
then $G$ is called a {\em md2 graph}. 
The transition matrix of a discrete-time quantum walk in a graph 
is closely related to the Ihara zeta function of a graph. 
We stare a relationship between the discrete-time quantum walk and the Ihara zeta function 
of a graph by Ren et al. \cite{RenETAL}.

\begin{thm}[Ren, Aleksic, Emms, Wilson and Hancock] 
Let ${\bf B} - {\bf J}_0 $ be the Perron-Frobenius operator (or the edge matrix) 
of a simple graph subject to the md2 constraint, where the edge matrix is defined in Section 3. 
Let ${\bf U}$ be the transition matrix of the discrete-time quantum walk on $G$. 
Then the ${\bf B} - {\bf J}_0 $ is the positive support of the transpose of ${\bf U} $, i.e., 
\[
{\bf B} - {\bf J}_0 =( {}^t {\bf U} )^+ . 
\]
\end{thm}

The structure theorem of the positive support of the square of the Grover matrix 
of a regular graph was obtained by Godsil and Guo \cite{GG2010}.

\begin{thm}[Godsil and Guo] 
Let $G$ be a connected $k$-regular graph with $m$ edges. 
Suppose that $k >2$. 
Then 
\[
( {\bf U}^2 )^+ =({\bf U}^+ )^2+{\bf I}_{2m} . 
\]
\end{thm}

Let $G$ be a connected graph and ${\bf U}$ the Grover transition matrix of $G$. 
A cycle $C=(v_0, e_1 , v_1, \cdots ,v_{n-1}, e_n , v_n )$ 
$(v_0 =v_n)$ is called {\em essential} if all vertices of $C$ 
except $v_0, v_n $ are distinct and $n \geq 2$. 
An essential cycle is the same as a dicycle in standard books on graph theory. 
Note that any essential cycle is a prime, reduced cycle, and 
any prime, reduced cycle is a union of disjoint essential cycles. 
The {\em girth} $g(G)$ of a graph $G$ is the minimum length of essential cycles in $G$. 

Higuchi, Konno, Sato, and Segawa \cite{HKSS} expressed $( {\bf U}^3)^+ $ by using ${\bf U}^+ $.

\begin{thm}[Higuchi, Konno, Sato, and Segawa]
Let $G$ be a $k$-regular graph. 
Suppose that $k>2$ and the girth $g(G)>4$.  
Then 
\[
({\bf U}^3)^+ =( {\bf U}^+ )^3 + {}^T {\bf U}^+ . 
\] 
\end{thm}

\section{The Ihara zeta function and the modified zeta function of a graph}

Let $G$ be a connected graph. Then a {\em path $P$ of length $n$} in $G$ is 
a sequence $P=(e_1, \ldots ,e_n )$ of $n$ arcs such that $e_i \in D(G)$, 
$t( e_i )=o( e_{i+1} )(1 \leq i \leq n-1)$, where indices are treated $mod \  n$. 
If $o(e_i)=v_{i-1} $ and $t(e_i )= v_i $ for $i=1, \ldots, n$, then we write 
$P=(v_0, e_1 , v_1, \cdots ,v_{n-1}, e_n , v_n )$.  
Set $ | P | =n$, $o(P)=o( e_1 )$ and $t(P)=t( e_n )$. 
Also, $P$ is called an {\em $(o(P),t(P))$-path}. 
Furthermore, $P$ is called an {\em $(e_1 , e_n )$-path}. 
We say that a path $P=(e_1, \ldots ,e_n )$ has a {\em backtracking} 
if $ e^{-1}_{i+1} =e_i $ for some $i (1 \leq i \leq n-1)$. 
A $(v, w)$-path is called a {\em $v$-cycle} 
(or {\em $v$-closed path}) if $v=w$. 
The {\em inverse cycle} of a cycle 
$C=( e_1, \ldots ,e_n )$ is the cycle 
$C^{-1} =( e^{-1}_n , \ldots ,e^{-1}_1 )$. 

We introduce an equivalence relation between cycles. 
Two cycles $C_1 =(e_1, \cdots ,e_m )$ and 
$C_2 =(f_1, \cdots ,f_m )$ are called {\em equivalent} if there exists 
$k$ such that $f_j =e_{j+k} $ for all $j$. 
The inverse cycle of $C$ is in general not equivalent to $C$. 
Let $[C]$ be the equivalence class which contains a cycle $C$. 
Let $B^r$ be the cycle obtained by going $r$ times around a cycle $B$. 
Such a cycle is called a {\em power} of $B$. 
A cycle $C$ is {\em reduced} if 
$C$ has no backtracking. 
Furthermore, a cycle $C$ is {\em prime} if it is not a power of 
a strictly smaller cycle. 
Note that each equivalence class of prime, reduced cycles of a graph $G$ 
corresponds to a unique conjugacy class of 
the fundamental group $ \pi {}_1 (G,v)$ of $G$ at a vertex $v$ of $G$. 

The {\em Ihara zeta function} of a graph $G$ is 
a function of $u \in {\bf C}$ with $\mid u \mid $ sufficiently small, 
defined by 
\[
{\bf Z} (G, u)= {\bf Z}_G (u)= \prod_{[C]} (1- u^{ \mid C \mid } )^{-1} ,
\]
where $[C]$ runs over all equivalence classes of prime, reduced cycles 
of $G$. 

Let $G$ be a connected graph with $n$ vertices and $m$ unoriented edges. 
Two $2m \times 2m$ matrices 
${\bf B} = {\bf B} (G)=( {\bf B}_{ef} )_{e,f \in D(G)} $ and 
${\bf J}_0 ={\bf J}_0 (G) =( {\bf J}_{ef} )_{e,f \in D(G)} $ 
are defined as follows: 
\[
{\bf B}_{ef} =\left\{
\begin{array}{ll}
1 & \mbox{if $t(e)=o(f)$, } \\
0 & \mbox{otherwise}
\end{array}
\right.
, 
{\bf J}_{ef} =\left\{
\begin{array}{ll}
1 & \mbox{if $f= e^{-1} $, } \\
0 & \mbox{otherwise.}
\end{array}
\right.
\]
Then the matrix ${\bf B} - {\bf J}_0 $ is called the {\em edge matrix} of $G$.

\begin{thm}[Hashimoto; Bass]
Let $G$ be a connected graph. 
Then the reciprocal of the zeta function of $G$ is given by 
\[
\begin{array}{rcl}
\  &   & {\bf Z} (G, u)^{-1} = \det ( {\bf I} -u ( {\bf B} - {\bf J}_0 )) \\
\  &   &                \\ 
\  & = & (1- u^2 )^{r-1} \det ( {\bf I} -u {\bf A} (G)+ u^2 ({\bf D} -{\bf I} ))
= \exp (\sum_{k \geq 1} \frac{N_k }{k} u^k ) , 
\end{array}  
\] 
where $r$ and ${\bf A} (G)$ are the Betti number and the adjacency matrix 
of $G$, respectively, and ${\bf D} ={\bf D} (G)=G( d_{ij} )$ is the diagonal matrix 
with $d_{ii} = \deg v_i $, $V(G)= \{ v_1 , \cdots , v_n \} $. 
Furthermore $N_k$ is the number of reduced cycles of length $k$ in $G$. 
\end{thm}

Next, we state a zeta function on the positive support of 
the Grover transition matrix of a graph. 
Let $G$ be a connected graph with n vertices and $m$ edges and 
${\bf U} = {\bf U} (G)$ the Grover transition matrix of $G$. 
By Theorems 2.3 and 3.1, we obtain the following result.

\begin{pro} 
Let $G$ be a connected graph with $m$ edges. 
Then  
\[
{\bf Z} (G, u)^{-1} = \det ( {\bf I}_{2m} -u {\bf U}^+ ) .  
\] 
\end{pro}

{\bf Proof}.  
\[
\det ( {\bf I}_{2m} -u {\bf U}^+ )^{-1} =\det ( {\bf I} -u ( {}^t {\bf B} - {}^t {\bf J}_0 ))
=\det ( {\bf I} -u ( {\bf B} - {\bf J}_0 ))^{-1} ={\bf Z} (G, u) . 
\] 
Q.E.D.

The Ihara zeta function of a graph is just a zeta function on the positive support of 
the Grover transition matrix of a graph. 

A zeta function of a graph related to the positive support of the square of 
the Grover transition matrix was defined by Higuch, Konno, Sato and Segawa  \cite{HKSS2014}.  
Let $G$ be a connected graph with n vertices and $m$ edges, $\delta (G)>2$ and 
${\bf U} = {\bf U} (G)$ the Grover transition matrix of $G$. 
Then the {\em modified zeta function} of $G$ is defined by  
\[
\tilde{{\bf Z}} (G,u)= \det ( {\bf I}_{2m} -u ({\bf U}^2)^+ )^{-1} .  
\]

Higuch, Konno, Sato and Segawa  \cite{HKSS2014} presented the Euler product, the exponential expression 
and another determinant expression for the modified zeta function of 
a graph.

\begin{thm}[Higuch, Konno, Sato and Segawa]  
Let $G$ be a connected graph with $n$ vertices and $m$ edges, 
let ${\bf U} = {\bf U} (G)$ be the Grover transition matrix of $G$. 
Suppose that $\delta (G)>2$.  
Then the modified zeta function of $G$ is given by 
\[
\tilde{{\bf Z}} (G,u)= \prod_{[C]} (1- u^{|C|} )^{-1} 
= \exp (\sum_{r \geq 1} \frac{N_r}{r} u^r ) 
\]
\[
=(1-2u )^{2(m-n)} \det ( {\bf I} -\sqrt{u(1-u)} {\bf A} (G)+u ({\bf D} -{\bf I} )) 
\det ( {\bf I} +\sqrt{u(1-u)} {\bf A} (G)+u ({\bf D} -{\bf I} )) . 
\]
where $[C]$ runs over all equivalence classes of prime 2-cycles in $G$, 
and $N_r $ is the number of 2-cycles of length $r$. 
\end{thm}

Thus, we propose the following problem:

\begin{prob}
Let $G$ be a connected graph with $m$ edges. 
Then, is there a zeta function $ \tilde{{\bf Z}} (G, u)$ of $G$ satisfying   
\[
{\bf Z}_3 (G, u)^{-1} = \det ( {\bf I}_{2m} -u ( {\bf U}^3 )^+ ) \   ?   
\]  
\end{prob}

Now, we consider a zeta function on the positive support of the cube of 
the Grover transition matrix of a graph.

\section{The structure theorem for the positive support of the cube 
of the Grover transition matrix of a general graph}

Now, we generalize Theorem 2.5 to a general graph.

\begin{thm}
Let $G$ be a connected graph.
Suppose that $\delta (G)>2$ and $g(G)>4$. 
The positive support $({\bf U}^3 ) ^+ $ is of the form 
\[
( {\bf U}^3 )^+ =( {\bf U}^+ )^3 + {}^T {\bf U}^+ . 
\]
\end{thm}

{\bf Proof}.  Now, let $G$ be a connected graph with $ \delta (G)>2$ and $g(G)>4$. 
Then we consider the structure of the positive support $( {\bf U}^3 )^+$ of the cube of 
the transition matrix ${\bf U} $.     
Since all nonzero elements of ${\bf B} $ and $ {}^T {\bf U} $ are in the same place, 
all nonzero elements of ${\bf B}^3 $ and $ {}^T {\bf U}^3 $ are in the same place. 
We treat ${\bf B}^3 $ and $ {}^T {\bf U}^3 $ in parallel. 

Let ${\bf T} = {\bf B} - {\bf J}_0 $ and ${\bf P} = {\bf J}_0 $. 
By Theorem 2.3, Then we have 
\[
{\bf B} = {\bf T} + {\bf P} \  and \  {}^T {\bf U}^+ = {\bf B} - {\bf P} . 
\]
Thus, we have 
\[
{\bf B}^3 = ( {\bf T} + {\bf P} )^3 
= {\bf T}^3 + {\bf T}^2 {\bf P} + {\bf T} {\bf P} {\bf T} + {\bf P} {\bf T}^2 + {\bf T} {\bf P}^2 
+ {\bf P}^2 {\bf T} + {\bf P}  {\bf T} {\bf P} + {\bf P}^3 . 
\]

But, the relation of arcs $e$ and $f$ of the nonzero $(e,f)$-array of $( {}^T {\bf U} )^3 $ are 
divided into the eight cases in Figure 1. 
In fact, the cases I, II, III, IV, V, VI, VII and VIII correspond to the matrices 
${\bf T}^3 $, ${\bf T}^2 {\bf P} $, $ {\bf T} {\bf P} {\bf T} $, ${\bf P} {\bf T}^2 $, 
$ {\bf T} {\bf P}^2 $, $ {\bf P}^2 {\bf T} $, ${\bf P}  {\bf T} {\bf P} $ and $ {\bf P}^3 $, 
respectively. 
In the case I, an $(e,f)$-path has no backtracking. 
For the cases II, III and IV, an $(e,f)$-path has exactly one backtracking. 
For the cases V, VI and VII, an $(e,f)$-path has exactly two backtrackings. 
In the case VIII, an $(e,f)$-path has exactly three backtrackings.  
The elements of ${}^T {\bf U}^3 $ corresponding to nonzero elements of 
${\bf T}^2 {\bf P} $, $ {\bf T} {\bf P} {\bf T} $, ${\bf P} {\bf T}^2 $ 
and $ {\bf P}^3 $ are negative. 
Furthermore, the elements of ${}^T {\bf U}^3 $ corresponding to nonzero elements of 
${\bf T}^3 $, $ {\bf T} {\bf P}^2 $, $ {\bf P}^2 {\bf T} $ and 
${\bf P}  {\bf T} {\bf P} $ are positive. 

If $t(e)=o(f)$, then nonzero $(e,f)$-arrays of $ {\bf T} {\bf P} {\bf T} $, $ {\bf T} {\bf P}^2 $ and 
$ {\bf P}^2 {\bf T} $ are overlapped. 
Then we have 
\[ 
\begin{array}{rcl}
( {\bf U}^3 )_{ef} & = & ( \frac{2}{d_{t(e)} } )^2 \sum_{o(g)=t(e), g \neq e,f}  ( \frac{2}{d_{t(e)}} -1) 
+ \frac{2}{d_{t(e)} } ( \frac{2}{d_{t(f)}} -1)( \frac{2}{d_{t(e)}} -1) \\
\  &  &    \\
\  & + & \frac{2}{d_{t(e)} } ( \frac{2}{d_{t(e)}} -1)( \frac{2}{d_{o(e)}} -1) . 
\end{array}
\]  
Since $\delta (G) \geq 3$, we have 
\[
d_{t(e)} \geq 3, \ i.e.,  \frac{1}{d_{t(e)}} \leq \frac{1}{3} 
\]
for any $e \in D(G)$. 
Thus, 
\[ 
\begin{array}{rcl}
( {\bf U}^3 )_{ef} & \leq  & ( \frac{2}{3} )^2 (d_{t(e)} -2)( \frac{2}{3} -1)+ \frac{2}{3} ( \frac{2}{3} -1)^2 \cdot 2 \\
\  &  &    \\
\  & = & - \frac{2}{9} ( \frac{2}{3} d_{t(e)} -2) \leq 0. 
\end{array}
\]
Therefore, all positive elements of $( {}^T {\bf U} )^3 $ and 
${\bf T}^3 + {\bf P} {\bf T} {\bf P}$ are in the same place. 
Therefore, 
\[
( {}^T {\bf U}^3 )^+ =( {\bf T}^3 + {\bf P}  {\bf T} {\bf P} )^+ . 
\]

Since $g(G)>4$, nonzero element of ${\bf T}^3 $ and ${\bf P}  {\bf T} {\bf P} $ are not overlapped. 
If a nonzero $(e,f)$-arrays of ${\bf T}^3 $ and ${\bf P} {\bf T} {\bf P}$ are overlapped, 
then there exists an essential cycle of length four from $e$ to $f$ in $G$, contradiction to $g(G)>4$. 

Furthermore, all nonzero elements of two matrices $ {\bf T}$ and ${\bf P}  {\bf T} {\bf P} $ are 1. 
In the case of $g(G)>4$, then all nonzero elements of the matrix $ {\bf T}^3 $ are 1. 
If an $(e,f)$-array of $ {\bf T}^3 $ is not less than 2, then there exist two distinct $(e,f)$-paths 
$P=(e,g,h,f)$ and $Q=(e,g_1 , h_1, f)$ in $G$. 
Then the cycle $(g,h,h^{-1}_1 , g^{-1}_1 )$ is an essential cycle of length four in $G$. 
This contradicts to the condition $g(G)>4$. 

Therefore, it follows that 
\[
( {\bf U}^3 )^+ =( {}^T {\bf T} )^3 + {\bf P} {}^T {\bf T} {\bf P} . 
\]
Since $ {}^T {\bf T} = {}^T {\bf B} - {\bf P} = {\bf U}^+ $, 
we have 
\[
( {\bf U}^3 )^+ =( {\bf U}^+ )^3 + {\bf P} {\bf U}^+ {\bf P} . 
\] 

But, by (5.3), 
\[
\begin{array}{rcl}
{\bf P} {\bf U}^+ {\bf P} & = & {\bf P} ( {}^T {\bf B} - {\bf P} ) {\bf P} 
={\bf P} ( {}^T {\bf L} {\bf K} - {\bf P} ) {\bf P} \\
\  &   &                \\ 
\  & = & {\bf P} {}^T {\bf L} {\bf K} {\bf P} - {\bf P}^3 
= {}^T {\bf K} {\bf L} - {\bf P} = {}^T {\bf U}^+ .  
\end{array}
\] 
Hence,   
\[
( {\bf U}^3 )^+ =( {\bf U}^+ )^3 + {}^T {\bf U}^+ . 
\]
Q.E.D.

\begin{center}
\begin{picture}(300,300)
\put(10,270){I}
\put(80,280){$e$}
\put(60,270){\circle{8}}
\put(65,275){\vector(1,0){30}}
\put(100,270){\circle{8}}
\put(105,275){\vector(1,0){30}}
\put(140,270){\circle{8}}
\put(145,275){\vector(1,0){30}}
\put(180,270){\circle{8}}
\put(185,275){\vector(1,0){30}}
\put(200,280){$f$}
\put(220,270){\circle{8}}
\put(10,240){II}
\put(80,250){$e$}
\put(60,240){\circle{8}}
\put(65,245){\vector(1,0){30}}
\put(100,240){\circle{8}}
\put(105,245){\vector(1,0){30}}
\put(140,240){\circle{8}}
\put(145,245){\vector(1,0){30}}
\put(180,240){\circle{8}}
\put(175,235){\vector(-1,0){30}}
\put(160,225){$f$}
\put(10,210){III}
\put(80,220){$e$}
\put(60,210){\circle{8}}
\put(65,215){\vector(1,0){30}}
\put(100,210){\circle{8}}
\put(105,215){\vector(1,0){30}}
\put(140,210){\circle{8}}
\put(135,205){\vector(-1,0){30}}
\put(100,205){\vector(0,-1){30}}
\put(100,170){\circle{8}}
\put(105,185){$f$}
\put(160,210){$\supset$}
\put(200,220){$e$}
\put(180,210){\circle{8}}
\put(185,215){\vector(1,0){30}}
\put(220,210){\circle{8}}
\put(215,205){\vector(-1,0){30}}
\put(200,185){$f$}
\put(225,215){\vector(1,0){30}}
\put(260,210){\circle{8}}
\put(255,205){\vector(-1,0){30}}
\put(10,150){IV}
\put(80,160){$f$}
\put(60,150){\circle{8}}
\put(95,155){\vector(-1,0){30}}
\put(100,150){\circle{8}}
\put(135,155){\vector(-1,0){30}}
\put(140,150){\circle{8}}
\put(160,160){$e$}
\put(145,155){\vector(1,0){30}}
\put(180,150){\circle{8}}
\put(175,145){\vector(-1,0){30}}
\put(10,120){V}
\put(80,130){$e$}
\put(60,120){\circle{8}}
\put(65,125){\vector(1,0){30}}
\put(100,120){\circle{8}}
\put(105,125){\vector(1,0){30}}
\put(140,120){\circle{8}}
\put(120,130){$f$}
\put(135,115){\vector(-1,0){30}}
\put(10,90){VI}
\put(80,100){$f$}
\put(60,90){\circle{8}}
\put(95,95){\vector(-1,0){30}}
\put(100,90){\circle{8}}
\put(135,95){\vector(-1,0){30}}
\put(140,90){\circle{8}}
\put(120,100){$e$}
\put(105,85){\vector(1,0){30}}
\put(10,60){VII}
\put(80,45){$f$}
\put(60,60){\circle{8}}
\put(65,55){\vector(1,0){30}}
\put(100,60){\circle{8}}
\put(95,65){\vector(-1,0){30}}
\put(105,65){\vector(1,0){30}}
\put(140,60){\circle{8}}
\put(120,70){$e$}
\put(135,55){\vector(-1,0){30}}
\put(10,30){VIII}
\put(80,37){$e$}
\put(60,30){\circle{8}}
\put(65,35){\vector(1,0){30}}
\put(100,30){\circle{8}}
\put(95,25){\vector(-1,0){30}}
\put(80,15){$f$} 
\end{picture}
\end{center}

\begin{center}
{\bf Figure 1} The nonzero $(e,f)$-array of $( {}^T {\bf U} )^3 $. 
\end{center}
    
\vspace{5mm}

\section{The zeta functions with respect to the positive support of the cube of 
the Grover transition matrix of a graph} 

In this section, we give a solution for Problem 3.4.
We define zeta function with respect to the positive support of the cube of the Grover matrix of a graph. 

At first, we introduce the notion of a new path and cycle in a graph. 
Let $G$ be a connected graph. 
A {\em 3-path $P$ of length $n$} in $G$ is a sequence 
$P=(e_1, \cdots ,e_{n} )$ of $n$ arcs such that $e_i \in D(G)$, and 
$( e_i , f, g, e_{i+1} ), f,g  \in D(G)$ is a reduced path of length three or 
($e_{i+1} , e_{i} ) \in D(G)$ for $i=1 , \ldots , n-1$. 
Here, $[ e_i , e_{i+1} ]=( e_i , f, g,  e_{i+1} )$ is called a {\em 3-arc}, and 
$[ e_i , e_{i+1} ]=( e_{i+1} , e_i )$is called {\em 3-backtracking}. 
Set $ \mid P \mid =n$, $o(P)=o(e_1 )$ and $t(P)=t(e_n )$. 
Also, $P$ is called an {\em $(o(P), t(P))$-3-path} or a {\em $(e_1 , e_{n} )$-3-path}. 
A $(v,w)$-path is called a {\em $v$-3-cycle} 
(or {\em $v$-3-closed path}) if $v=w$. 
The {\em inverse cycle} of a 3-cycle 
$C=( e_1, \cdots ,e_n )$ is the broad cycle  
$C^{-1} =( e^{-1}_n, \cdots ,e^{-1}_1 )$. 

We introduce an equivalence relation between 3-cycles. 
Two 3-cycles $C_1 =(e_1, \cdots ,e_m )$ and 
$C_2 =(f_1, \cdots ,f_m )$ are called {\em equivalent} if 
$f_j =e_{j+k} $ for all $j$. 
The inverse 3-cycle of $C$ is in general not equivalent to $C$. 
Let $[C]$ be the equivalence class which contains a 3-cycle $C$. 
Let $B^r$ be the 3-cycle obtained by going $r$ times around a broad cycle $B$. 
A 3-cycle $C$ is {\em prime} if it is not a multiple of 
a strictly smaller 3-cycle.

Next, we define a zeta function of a graph related to the positive support of the cube of 
the Grover transition matrix. 
Let $G$ be a connected graph with n vertices and $m$ edges, $\delta (G)>2$ and 
${\bf U} = {\bf U} (G)$ the Grover transition matrix of $G$. 
Suppose that $g(G)>4$. 
Then the {\em 3-modified zeta function} of $G$ is defined by  
\[
{\bf Z}_3 (G,u)= \det ( {\bf I}_{2m} -u ({\bf U}^3)^+ )^{-1} .  
\]

We present the Euler product, the exponential expression 
and another determinant expression for the 3-modified zeta function of 
a graph.

\begin{thm}
Let $G$ be a connected graph with $n$ vertices and $m$ edges, 
let ${\bf U} = {\bf U} (G)$ be the Grover transition matrix of $G$. 
Suppose that $ \delta (G)>2$ and $g(G)>4$.  
Then the 3-modified zeta function of $G$ is given by 
\[
\begin{array}{rcl}
\  &  & {\bf Z}_3 (G,u)^{-1} = \prod_{[C]} (1- u^{|C|} ) 
= \exp (- \sum_{r \geq 1} \frac{N_r}{r} u^r ) \\
\  &   &                \\ 
\  & = & (1-4 u^2 )^{m-2n} \det ( {\bf I}_{n} +u( {\bf D} -2 {\bf I}_n ) {\bf A} -2u^2 ( {\bf A}^2 -2 {\bf D} +2 {\bf I}_n ) \\  
\  &   &                \\ 
\  & + & u^2 ( ({\bf D} - {\bf I}_n ) {\bf A}^2 - {\bf D}^2 +2 {\bf D} -2u ( {\bf A}^3 -({\bf D} -2 {\bf I}_n ) {\bf A} - {\bf A} {\bf D} )) \\ 
\  &   &                \\ 
\  & \times & ( {\bf I}_n -u( {\bf A}^3 -({\bf D} -2 {\bf I}_n ) {\bf A} - {\bf A} {\bf D} )
+2u^2 (( {\bf D} -{\bf I}_n ) {\bf A}^2 - {\bf D}^2 +2 {\bf D} -2 {\bf I}_n ))^{-1} \\ 
\  &   &                \\ 
\  & \times & ( {\bf A}^2 -2u( {\bf D} -2 {\bf I}_n ) {\bf A} -2 {\bf D} )) \\
\  &   &                \\ 
\  & \times & \det ( {\bf I}_{n} -u( {\bf A}^3 -({\bf D} -2 {\bf I}_n ) {\bf A} - {\bf A} {\bf D} )
+2u^2 (( {\bf D} -{\bf I}_n ) {\bf A}^2 - {\bf D}^2 +2 {\bf D} -2 {\bf I}_n )) , 
\end{array}
\]
where $[C]$ runs over all equivalence classes of prime 3-cycles in $G$, 
and $N_r $ is the number of 3-cycles of length $r$. 
\end{thm}

{\bf Proof}.  We use  the method of Smilansky \cite{Sm}. 
By Theorem 4.1, we have 
\[
{\bf Z}_3 (G,u)^{-1} = \det ( {\bf I}_{2m} -u (({\bf U}^+)^3 + {}^T {\bf U}^+ ))
= \det ( {\bf I}_{2m} -u (({\bf U}^+)^3 + {}^T {\bf U}^+ )) .  
\] 

At first, 
\[
\log \det ( {\bf I}_{2m} -u (({\bf U}^+)^3 + {}^T {\bf U}^+ ))
= {\rm Tr} \log ( {\bf I}_{2m} -u (({\bf U}^+)^3 + {}^T {\bf U}^+ ))
=- \sum^{\infty}_{n=1} \frac{1}{n} {\rm Tr} [({\bf U}^+)^3 + {}^T {\bf U}^+ )^n ) u^n ] . 
\]

But, we have 
\[
{\rm Tr} [({\bf U}^+)^3 + {}^T {\bf U}^+ )^n u^n ] 
\]
\[
= u^n |\{ C_n =(e_1, \cdots ,e_{n} , e_1 ) \mid [e_i ,e_{i+1} ]: 3-arc \  or \  3-backtracking \} | . 
\]
For a 3-cycle $C_n$ with length $n$, the exists a prime 3-cycle $\tilde{C}_p $ such that 
$C_n = \tilde{C}_p $ and $n=pk$. 
Thus, we have 
\[
{\rm Tr} [(({\bf U}^+)^3 + {}^T {\bf U}^+ )^n u^n ] 
= \sum_{[\tilde{C}_p]} p u^{pk} ,  
\]
where $[\tilde{C}_p]$ runs over all equivalence classes of prime 3-cycles in $G$, 
Therefore, it follows that 
\[
\log \det ( {\bf I}_{2m} -u (({\bf U}^+)^3 + {}^T {\bf U}^+ ))
=- \sum_{[\tilde{C}_p]} \sum^{\infty}_{k=1} \frac{1}{kp} p u^{pk} 
= \sum_{[\tilde{C}_p]} \log (1- u^{|\tilde{C}_p |} ) . 
\]
Hence, 
\[
{\bf Z}_3 (G,u)= \prod_{{\tilde{C}_p }} (1- u^{|\tilde{C}_p |} )^{-1} . 
\]

Next, we give the exponential expression of the 3-modified zeta function. 
By the definition of $N_k $, we have 
\[
N_k = {\rm Tr} [({\bf U}^+)^3 + {}^T {\bf U}^+ )^k ]. 
\]
Then we have 
\[
\log {\bf Z}_3 (G,u)^{-1} =- \sum^{\infty}_{k=1} \frac{N_k}{k} u^k .
\]
Thus, 
\[
{\bf Z}_3 (G,u)= \exp (\sum_{k \geq 1} \frac{N_k}{k} u^k ) . 
\]

Let $V(G)= \{ v_1 , \ldots , v_n \}$ and 
$D(G)= \{ e_1 , \ldots , e_{m} $, $e^{-1}_1 , \ldots , e^{-1}_{m} \} $. 
Arrange arcs of $G$ as follows:
$e_1 , e^{-1}_1 , \ldots , e_{m}, e^{-1}_m $. 
Furthermore, arrange vertices of $G$ as follows:
$ v_1 , \ldots , v_n $.

Now, we define two $n \times 2m$ matrices 
${\bf K} =( {\bf K}_{ve} ) {}_{e \in D(G); v \in V(G)} $ and 
${\bf L} =( {\bf L}_{ve} )_{e \in D(G); v \in V(G)} $ as follows: 
\[
{\bf K}_{ve} :=\left\{
\begin{array}{ll}
1 & \mbox{if $t(e)=v$, } \\
0 & \mbox{otherwise. } 
\end{array}
\right. 
,
{\bf L}_{ve} :=\left\{
\begin{array}{ll}
1 & \mbox{if $o(e)=v$, } \\
0 & \mbox{otherwise. } 
\end{array}
\right.
\]
Here we consider two matrices ${\bf K}$ and  ${\bf L}$ under the above order. 
Then we have 
\begin{equation}
{\bf L} {}^t {\bf K} = {\bf K} {}^t {\bf L} ={\bf A} (G) , 
\end{equation} 
\begin{equation}
{\bf L} {}^t {\bf L} = {\bf K} {}^t {\bf K} ={\bf D} ,  
\end{equation} 
\begin{equation}
{}^t {\bf K} {\bf L} = {\bf B} (G) \  and \  
{}^t {\bf L} {\bf K} = {}^T {\bf B} (G) . 
\end{equation}
Furthermore, 
\begin{equation}
{\bf L} = {\bf K} {\bf J}_0 \  and \  {\bf K} = {\bf L} {\bf J}_0 , 
\end{equation} 
\begin{equation}
{}^T {\bf L} = {\bf J}_0 {}^T {\bf K} \  and \  {}^T {\bf K} = {\bf J}_0 {}^T {\bf L} .  
\end{equation}
Note that, if $\delta (G) \geq 2$, then  
\begin{equation}
{\bf U}^+ ={}^T {\bf B} - {\bf J}_0 . 
\end{equation}

Now, set 
\[
{\bf F} = {}^T {\bf B} (G)= {}^t {\bf L} {\bf K} \  and \ {\bf J} = {\bf J}_0 . 
\]
Then we have 
\[
\begin{array}{rcl}
\  &  & ( {\bf U}^+)^3 =( {}^T {\bf B} - {\bf J}_0 )^3 =( {\bf F} - {\bf J} )^3 \\
\  &   &                \\ 
\  & = & {\bf F}^3 - {\bf F} {\bf J} {\bf F} - {\bf J} {\bf F}^2 +2 {\bf F} 
- {\bf F} {}^2 {\bf J} + {\bf J} {\bf F} {\bf J} - {\bf J} \\ 
\  &   &                \\ 
\  & = &  {}^t {\bf L} {\bf K} {}^t {\bf L} {\bf K}{}^t {\bf L} {\bf K} - {}^t {\bf L} {\bf K} {\bf J} {}^t {\bf L} {\bf K} 
- {\bf J} {}^t {\bf L} {\bf K} {}^t {\bf L} {\bf K} 
+2  {}^t {\bf L} {\bf K} - {}^t {\bf L} {\bf K}  {}^t {\bf L} {\bf K} {\bf J} 
+ {\bf J} {}^t {\bf L} {\bf K} {\bf J} - {\bf J} \\ 
\  &   &                \\ 
\  & = &  {}^t {\bf L} {\bf A} {}^2 {\bf K} - {}^t {\bf L} {\bf K} {}^t {\bf K} {\bf K} 
- {}^t {\bf K} {\bf A} {\bf K} 
+2  {}^t {\bf L} {\bf K} -  {}^t {\bf L} {\bf A} {\bf L} + {}^t {\bf K} {\bf L} - {\bf J} \\ 
\  &   &                \\ 
\  & = & {}^t {\bf L} {\bf A} {}^2 {\bf K} -k {}^t {\bf L} {\bf K} - {}^t {\bf K} {\bf A} {\bf K} 
+2 {}^t {\bf L} {\bf K} -  {}^t {\bf L} {\bf A} {\bf L} + {}^t {\bf K} {\bf L} - {\bf J} . 
\end{array}
\]    
Furthermore, ewe have 
\[
{}^t {\bf U}^+ = {}^T {\bf K} {\bf L} - {\bf J} . 
\]
Thus, 
\[
( {\bf U}^+)^3 + {}^T {\bf U}^+ = 
{}^t {\bf L} ( {\bf A} {}^2 -(k-2) {\bf I}_n ) {\bf K} - {}^t {\bf K} {\bf A} {\bf K} 
-  {}^t {\bf L} {\bf A} {\bf L} +2 {}^t {\bf K} {\bf L} -2 {\bf J}  . 
\]    
Therefore, it follows that 
\begin{equation}
\begin{array}{rcl}
\  &  & \det ( {\bf I}_{2m} -u ( {\bf U}^3)^+ )\\ 
\  &   &                \\ 
\  & = & \det ( {\bf I}_{2m} -u( {}^t {\bf L} ( {\bf A} {}^2 -(k-2) {\bf I}_n ) {\bf K} - {}^t {\bf K} {\bf A} {\bf K} 
-  {}^t {\bf L} {\bf A} {\bf J} +2 {}^t {\bf K} {\bf L} -2 {\bf J} )) \\
\  &   &                \\ 
\  & = & \det ( {\bf I}_{2m} +2u {\bf J}-u( {}^t {\bf L} ( {\bf A} {}^2 -(k-2) {\bf I}_n ) {\bf K} - {}^t {\bf K} {\bf A} {\bf K} 
-  {}^t {\bf L} {\bf A} {\bf L} +2 {}^t {\bf K} {\bf L} )) \\
\  &   &                \\ 
\  & = & \det ( {\bf I}_{2m} -u( {}^t {\bf L} ( {\bf A} {}^2 -(k-2) {\bf I}_n ) {\bf K} - {}^t {\bf K} {\bf A} {\bf K} 
-  {}^t {\bf L} {\bf A} {\bf L} +2 {}^t {\bf K} {\bf L} ) \\
\  &   &                \\ 
\  & \times & ( {\bf I}_{2m} +2u {\bf J} )^{-1} ) \det ( {\bf I}_{2m} +2u {\bf J} ) . 
\end{array}
\end{equation}

But, we have 
\[
\det ( {\bf I}_{2m} +2u {\bf J} )=(1-4 u^2 )^m . 
\]
Furthermore, we have 
\[
( {\bf I}_{2m} +2u {\bf J} )^{-1} = \frac{1}{1-4 u^2 }  
\left[ 
\begin{array}{cccc}
1  &  -2u &   & 0 \\
-2u & 1 &   &  \\ 
0 &    & \ddots &  
\end{array} 
\right] 
= \frac{1}{1-4 u^2 }( {\bf I}_{2m} -2u {\bf J} ) .  
\]
Thus, 
\[
\begin{array}{rcl}
\  &  & (( {}^t {\bf L} ( {\bf A} {}^2 -(k-2) {\bf I}_n ) {\bf K} - {}^t {\bf K} {\bf A} {\bf K} 
-  {}^t {\bf L} {\bf A} {\bf L} +2 {}^t {\bf K} {\bf L} )( {\bf I}_{2m} +2u {\bf J} )^{-1} \\ 
\  &   &                \\ 
\  & = & \frac{1}{1-4 u^2 } 
( {}^t {\bf L} ( {\bf A} {}^2 -(k-2) {\bf I}_n ) {\bf K} - {}^t {\bf K} {\bf A} {\bf K} 
-  {}^t {\bf L} {\bf A} {\bf L} +2 {}^t {\bf K} {\bf L} ) \\
\  &   &                \\ 
\  & - & 2u( {}^t {\bf L} ( {\bf A} {}^2 -(k-2) {\bf I}_n ) {\bf K} {\bf J} - {}^t {\bf K} {\bf A} {\bf K} {\bf J} 
-  {}^t {\bf L} {\bf A} {\bf L} {\bf J} -2 {}^t {\bf K} {\bf L}{\bf J} )) \\
\  &   &                \\ 
\  & = & \frac{1}{1-4 u^2 } 
( \{ {}^t {\bf L} ( {\bf A} {}^2 -(k-2) {\bf I}_n )- {}^t {\bf K} {\bf A} \} ( {\bf K} -2u {\bf L} )  
- ( {}^t {\bf L} {\bf A} -2 {}^t {\bf K} )( {\bf L} -2u {\bf K} )) . 
\end{array}
\]

Next, let 
\[
x= \frac{u}{1-4 u^2} , \ {\bf P} ={}^t {\bf L} ( {\bf A} {}^2 -(k-2) {\bf I}_n )- {}^t {\bf K} {\bf A} , 
\]
\[
{\bf Q} = {\bf K} -2u {\bf L} , {\bf R} = {}^t {\bf L} {\bf A} -2 {}^t {\bf K}, 
{\bf T} = {\bf L} -2u {\bf K} . 
\]
Then we have 
\[
\begin{array}{rcl}
\  &  &  \det ( {\bf I}_{2m} -u ( {\bf U}^3)^+ )\\ 
 &   &                \\ 
\  & = & \det ({\bf I}_{2m} -x( {\bf P} {\bf Q} - {\bf R} {\bf T}))(1-4 u^2 )^{m} \\ 
\  &   &                \\ 
\  & = & \det ({\bf I}_{2m} -x {\bf P} {\bf Q} +x {\bf R} {\bf T}))(1-4 u^2 )^{m} \\ 
 &   &                \\ 
\  & = & \det ({\bf I}_{2m} +x {\bf R} {\bf T} ({\bf I}_{2m} -x {\bf P} {\bf Q})^{-1} ) 
\det ({\bf I}_{2m} -x {\bf P} {\bf Q})(1-4 u^2 )^{m} . 
\end{array}
\]

Furthermore, if ${\bf A}$ and ${\bf B}$ are a $m \times n $ and $n \times m$ 
matrices, respectively, then we have 
\[
\det ( {\bf I}_m - {\bf A} {\bf B} )= 
\det ( {\bf I}_n - {\bf B} {\bf A} ) . 
\]
Thus, we have 
\[
\begin{array}{rcl}
\  &  &  \det ( {\bf I}_{2m} -u ( {\bf U}^3)^+ )\\ 
\  &   &                \\ 
\  & = & \det ({\bf I}_{n} +x {\bf T} ({\bf I}_{2m} -x {\bf P} {\bf Q})^{-1}{\bf R} ) 
\det ({\bf I}_{n} -x {\bf Q}{\bf P} )(1-4 u^2 )^{m} \\ 
\  &   &                \\ 
\  & = & \det ({\bf I}_{n} +x {\bf T} ({\bf I}_{2m} +x {\bf P} {\bf Q} +x^2 {\bf P} {\bf Q} {\bf P} {\bf Q} 
+ x^3 {\bf P} {\bf Q} {\bf P} {\bf Q} {\bf P} {\bf Q} + \cdots ) {\bf R} )  \\
\  &   &                \\ 
\  & \times & 
\det ({\bf I}_{n} -x {\bf Q}{\bf P} )(1-4 u^2 )^{m} \\ 
\  &   &                \\ 
\  & = & \det ({\bf I}_{n} +x {\bf T} {\bf R} +x^2 {\bf T} {\bf P} {\bf Q} {\bf R} +x^3 {\bf T} {\bf P} {\bf Q} {\bf P} {\bf Q} {\bf R} 
+ x^4 {\bf T} {\bf P} {\bf Q} {\bf P} {\bf Q} {\bf P} {\bf Q} {\bf R} + \cdots ) \\ 
\  &   &                \\ 
\  & \times & \det ({\bf I}_{n} -x {\bf Q}{\bf P} )(1-4 u^2 )^{m} \\ 
\  &   &                \\ 
\  & = & \det ({\bf I}_{n} +x {\bf T} {\bf R} +x^2 {\bf T} {\bf P} ( {\bf I}_n +x {\bf Q} {\bf P} +x^2 ( {\bf Q} {\bf P} )^2 + \cdots ) 
{\bf Q} {\bf R} ) \\ 
\  &   &                \\ 
\  & \times & \det ({\bf I}_{n} -x {\bf Q}{\bf P} )(1-4 u^2 )^{m} \\ 
\  &   &                \\ 
\  & = & \det ({\bf I}_{n} +x {\bf T} {\bf R} +x^2 {\bf T} {\bf P} ( {\bf I}_n -x {\bf Q} {\bf P} )^{-1} {\bf Q} {\bf R} ) 
\det ({\bf I}_{n} -x {\bf Q}{\bf P} )(1-4 u^2 )^{m} . 
\end{array}
\]
Since $x= \frac{u}{1 -4 u^2} $. 
Then we have 
\[
\begin{array}{rcl}
\  &  &  \det ( {\bf I}_{2m} -u ( {\bf U}^3)^+ )\\ 
\  &   &                \\ 
\  & = & \det ({\bf I}_{n} +x {\bf T} {\bf R} +x^2 {\bf T} {\bf P} ( {\bf I}_n -x {\bf Q} {\bf P} )^{-1} {\bf Q} {\bf R} ) 
\det ({\bf I}_{n} -x {\bf Q}{\bf P} )(1-4 u^2 )^{m} \\
\  &   &                \\ 
\  & = & \det ({\bf I}_{n} + \frac{u}{1-4 u^2} {\bf T} {\bf R} +( \frac{u}{1-4 u^2} )^2 {\bf T} {\bf P} 
( {\bf I}_n - \frac{u}{1-4 u^2} {\bf Q} {\bf P} )^{-1} {\bf Q} {\bf R} ) \\
\  &   &                \\ 
\  & \times & \det ({\bf I}_{n} - \frac{u}{1-u^2} {\bf Q}{\bf P} )(1-4 u^2 )^{m} \\
\  &   &                \\ 
\  & = & \det ((1 -4 u^2 ) {\bf I}_{n} +u {\bf T} {\bf R} +u^2 {\bf T} {\bf P} 
((1 -4 u^2 ) {\bf I}_n -u {\bf Q} {\bf P} )^{-1} {\bf Q} {\bf R} ) \\
\  &   &                \\ 
\  & \times & \det ((1 -4 u^2 ) {\bf I}_{n} -u {\bf Q}{\bf P} )(1-4 u^2 )^{m-2n} . 
\end{array}
\]

But, we have 
\begin{equation}
\begin{array}{rcl}
\  &  & {\bf T} {\bf R} =({\bf L} -2u {\bf K} )( {}^t {\bf L} {\bf A} -2 {}^t {\bf K} ) \\
 &   &                \\ 
\  & = & {\bf L} {}^t {\bf L} {\bf A} -2 {\bf L} {}^T {\bf K} -2u {\bf K} {}^t {\bf L} {\bf A} +4u {\bf K} {}^t {\bf K} \\ 
 &   &                \\ 
\  & = & {\bf D} {\bf A} -2 {\bf A} -2u {\bf A}^2 +4u {\bf D} =-2u {\bf A}^2 +( {\bf D} -2 {\bf I}_n ) {\bf A} +4u {\bf D} ,  
\end{array}
\end{equation} 
\begin{equation}
\begin{array}{rcl}
\  &  & {\bf T} {\bf P} =({\bf L} -2u {\bf K} )( {}^t {\bf L} ( {\bf A} {}^2 -(k-2) {\bf I}_n )- {}^t {\bf K} {\bf A} ) \\
 &   &                \\ 
\  & = & {\bf L} {}^t {\bf L} ({\bf A}^2 -{\bf D} +2 {\bf I}_n )- {\bf L} {}^t {\bf K} {\bf A} 
-2u {\bf K} {}^t {\bf L} ({\bf A}^2 -{\bf D} +2 {\bf I}_n )+2u {\bf K} {}^T {\bf K} {\bf A} \\
 &   &                \\ 
\  & = & {\bf D} {\bf A}^2 - {\bf D}^2 +2 {\bf D} - {\bf A}^2 -2u {\bf A}^3 +2u {\bf A} {\bf D} 
-4u {\bf A} +2u {\bf D} {\bf A} \\
\  &   &                \\ 
\  & = & -2u {\bf A}^3 +({\bf D} - {\bf I}_n ) {\bf A}^2 +2u({\bf D} -2 {\bf I}_n ) {\bf A} - {\bf D}^2 +2( {\bf I}_n +u {\bf A} ) {\bf D} ,  
\end{array}
\end{equation} 
\begin{equation}
\begin{array}{rcl}
\  &  & {\bf Q} {\bf P} =({\bf K} -2u {\bf L} )( {}^t {\bf L} ( {\bf A} {}^2 -{\bf D} +2 {\bf I}_n )- {}^t {\bf K} {\bf A} ) \\
 &   &                \\ 
\  & = & {\bf K} {}^t {\bf L} ({\bf A}^2 -{\bf D} +2 {\bf I}_n )- {\bf K} {}^t {\bf K} {\bf A} 
-2u {\bf L} {}^t {\bf L} ({\bf A}^2 -{\bf D} +2 {\bf I}_n )+2u {\bf L} {}^T {\bf K} {\bf A} \\
\  &   &                \\ 
\  & = & {\bf A}^3 - {\bf A} {\bf D} +2 {\bf A} - {\bf D} {\bf A} -2u {\bf D} {\bf A}^2 +2u {\bf D}^2 -4u {\bf D} +2u {\bf A}^2  \\
\  &   &                \\ 
\  & = & {\bf A}^3 -2u({\bf D} - {\bf I}_n ) {\bf A}^2 -({\bf D} -2 {\bf I}_n ) {\bf A} +2u {\bf D}^2 -( {\bf A} +4u {\bf I}_n ) {\bf D} ,   
\end{array}
\end{equation}
\begin{equation}
\begin{array}{rcl}
\  &  & {\bf Q} {\bf R} =({\bf K} -2u {\bf L} )( {}^t {\bf L} {\bf A} -2 {}^t {\bf K} ) \\
\  &   &                \\ 
\  & = & {\bf K} {}^t {\bf L} {\bf A} -2 {\bf K} {}^T {\bf K} -2u {\bf L} {}^t {\bf L} {\bf A} +4u {\bf L} {}^t {\bf K} \\ 
\  &   &                \\ 
\  & = & {\bf A}^2 -2 {\bf D} -2u {\bf D} {\bf A} +4u {\bf A} 
= {\bf A}^2 -2u( {\bf D} -2 {\bf I}_n ) {\bf A} -2 {\bf D} . 
\end{array}
\end{equation}
Therefore, it follows that 
\[
\begin{array}{rcl}
\  &  &  \det ( {\bf I}_{2m} -u ( {\bf U}^3)^+ ) \\ 
\  &   &                \\ 
\  & = & \det ((1 -4 u^2 ) {\bf I}_{n} +u( -2u {\bf A}^2 +( {\bf D} -2 {\bf I}_n ) {\bf A} +4u {\bf D} ) \\  
\  &   &                \\ 
\  & + & u^2 ( -2u {\bf A}^3 +({\bf D} - {\bf I}_n ) {\bf A}^2 +2u({\bf D} -2 {\bf I}_n ) {\bf A} - {\bf D}^2 +2( {\bf I}_n +u {\bf A} ) {\bf D} ) \\ 
\  &   &                \\ 
\  & \times & ((1 -4 u^2 ) {\bf I}_n -u( {\bf A}^3 -2u({\bf D} - {\bf I}_n ) {\bf A}^2 -({\bf D} -2 {\bf I}_n ) {\bf A} 
+2u {\bf D}^2 -( {\bf A} +4u {\bf I}_n ) {\bf D} ) )^{-1}  \\ 
\  &   &                \\ 
\  & \times & ( {\bf A}^2 -2u( {\bf D} -2 {\bf I}_n ) {\bf A} -2 {\bf D} )) \\
\  &   &                \\ 
\  & \times & \det ((1 -4 u^2 ) {\bf I}_{n} -u( {\bf A}^3 -2u({\bf D} - {\bf I}_n ) {\bf A}^2 -({\bf D} -2 {\bf I}_n ) {\bf A} 
+2u {\bf D}^2 -( {\bf A} +4u {\bf I}_n ) {\bf D} )) \\
\  &   &                \\ 
\  & \times & (1-4 u^2 )^{m-2n} \\  
\  &   &                \\ 
\  & = & \det ( {\bf I}_{n} +u( {\bf D} -2 {\bf I}_n ) {\bf A} -2u^2 ( {\bf A}^2 -2 {\bf D} +2 {\bf I}_n ) \\  
\  &   &                \\ 
\  & + & u^2 ( ({\bf D} - {\bf I}_n ) {\bf A}^2 - {\bf D}^2 +2 {\bf D} -2u ( {\bf A}^3 -({\bf D} -2 {\bf I}_n ) {\bf A} - {\bf A} {\bf D} )) \\ 
\  &   &                \\ 
\  & \times & ( {\bf I}_n -u( {\bf A}^3 -({\bf D} -2 {\bf I}_n ) {\bf A} - {\bf A} {\bf D} )
+2u^2 (( {\bf D} -{\bf I}_n ) {\bf A}^2 - {\bf D}^2 +2 {\bf D} -2 {\bf I}_n ))^{-1} \\ 
\  &   &                \\ 
\  & \times & ( {\bf A}^2 -2u( {\bf D} -2 {\bf I}_n ) {\bf A} -2 {\bf D} )) \\
\  &   &                \\ 
\  & \times & \det ( {\bf I}_{n} -u( {\bf A}^3 -({\bf D} -2 {\bf I}_n ) {\bf A} - {\bf A} {\bf D} )
+2u^2 (( {\bf D} -{\bf I}_n ) {\bf A}^2 - {\bf D}^2 +2 {\bf D} -2 {\bf I}_n )) \\
\  &   &                \\ 
\  & \times & (1-4 u^2 )^{m-2n} .   
\end{array}
\]
Q.E.D.

In the case that $G$ is regular, we obtain the following result from Theorem 5.1.

\begin{thm}
Let $G$ be a connected $k$-regular graph with $n$ vertices 
and $m$ unoriented edges, $k>2$ and $g(G)>4$.  
Then 
\[
\begin{array}{rcl}
{\bf Z}_3 (G,u)^{-1} & = & (1-4 u^2 )^{m-n}  
\det ( {\bf I}_{n} -u ( {\bf A}^3 -(3k-4) {\bf A} ) \\
\  &   &                \\ 
\  & + & 
u^2 ( {\bf A}^4 - k^2 {\bf A}^2 +2(k-1)(k^2 -2k+2) {\bf I}_n )) .
\end{array}
\]    
\end{thm}

{\bf Proof}.  Let $G$ be $k$-regular. 
Then we have 
\[
{\bf D} =k {\bf I}_n . 
\]
Thus, 
\[
\begin{array}{rcl}
\  &  &  \det ( {\bf I}_{2m} -u ( {\bf U}^3)^+ ) \\ 
\  &   &                \\ 
\  & = & \det ( {\bf I}_{n} +u(k-2) {\bf A} -2u^2 ( {\bf A}^2 -2(k-1) {\bf I}_n ) \\  
\  &   &                \\ 
\  & + & u^2 ( (k-1) {\bf A}^2 -k(k-2) {\bf I}_n -2u ( {\bf A}^3 -(k-2) {\bf A} -k {\bf A} )) \\ 
\  &   &                \\ 
\  & \times & ( {\bf I}_n -u( {\bf A}^3 -(k-2) {\bf A} -k {\bf A} )
+2u^2 ((k-1) {\bf A}^2 -(k^2 -2k+2) {\bf I}_n ))^{-1} \\ 
\  &   &                \\ 
\  & \times & ( {\bf A}^2 -2u(k-2) {\bf A} -2k {\bf I}_n )) \\
\  &   &                \\ 
\  & \times & \det ( {\bf I}_{n} -u( {\bf A}^3 -(k-2) {\bf A} -k {\bf A} )
+2u^2 ((k-1) {\bf A}^2 -(k^2 -2k+2) {\bf I}_n )) \\
\  &   &                \\ 
\  & \times & (1-4 u^2 )^{m-2n}  .   
\end{array}
\]

Since the polynomials of ${\bf A} $ are commutative, we have 
\[ 
\begin{array}{rcl}
\  &  &  \det ( {\bf I}_{2m} -u ( {\bf U}^3)^+ ) \\ 
\  &   &                \\ 
\  & = & \det ( {\bf I}_{n} +u(k-2) {\bf A} -2u^2 ( {\bf A}^2 -2(k-1) {\bf I}_n ) ) \\  
\  &   &                \\ 
\  & \times & ( {\bf I}_{n} -u( {\bf A}^3 -(k-2) {\bf A} -k {\bf A} )
+2u^2 ((k-1) {\bf A}^2 -(k^2 -2k+2) {\bf I}_n )) \\
\  &   &                \\ 
\  & + & u^2 ( (k-1) {\bf A}^2 -k(k-2) {\bf I}_n -2u ( {\bf A}^3 -(k-2) {\bf A} -k {\bf A} )) \\ 
\  &   &                \\ 
\  & \times & ( {\bf A}^2 -2u(k-2) {\bf A} -2k {\bf I}_n )) (1-4 u^2 )^{m-2n} \\ 
\  &   &                \\ 
\  & = & \det ( {\bf I}_{n} -u( {\bf A}^3 -(3k-4) {\bf A} )+u^2 ( {\bf A}^4 - k^2 {\bf A}^2 +(2k^3 -6k^2 +8k-8) {\bf I}_n ) \\  
\  &   &                \\ 
\  & + & 4u^3 ( {\bf A}^3 -(3k-4) {\bf A} )+4u^4 ( - {\bf A}^4 + k^2 {\bf A}^2 )+2(k-1)(k^2-2k+2) {\bf I}_n )) \\ 
\  &   &                \\ 
\  & \times & (1-4 u^2 )^{m-2n} \\
\  &   &                \\ 
\  & = & \det ((1-4 u^2 )( {\bf I}_{n} -u( {\bf A}^3 -(3k-4) {\bf A} )+u^2 ( {\bf A}^4 - k^2 {\bf A}^2 +2(k-1)(k^2-2k+2) {\bf I}_n )) \\  
\  &   &                \\ 
\  & \times & (1-4 u^2 )^{m-2n} \\
\  &   &                \\ 
\  & = & \det ( {\bf I}_{n} -u( {\bf A}^3 -(3k-4) {\bf A} )+u^2 ( {\bf A}^4 - k^2 {\bf A}^2 +2(k-1)(k^2-2k+2) {\bf I}_n ))) \\  
\  &   &                \\ 
\  & \times & (1-4 u^2 )^{m-n} . 
\end{array}
\]
Q.E.D.

Thus, 

\begin{cor} 
Let $G$ be a connected $k$-regular graph with $n$ vertices 
and $m$ unoriented edges, $k>2$ and $g(G)>4$.  
Then the characteristic polynomial of $( {\bf U}^3)^+ $ is 
\[
\det ( \lambda {\bf I}_{2m} - ( {\bf U}^3)^+ )
=( \lambda {}^2 -4 )^{m-n} \det (  \lambda {}^2 {\bf I}_{n} - \lambda ( {\bf A}^3 -(3k-4) {\bf A} )
+( {\bf A}^4 - k^2 {\bf A}^2 +2(k-1)(k^2 -2k+2) {\bf I}_{n} )) .   
\] 
\end{cor}

Furthermore, we obtain the spectra for  $( {\bf U}^3)^+ $.

\begin{cor} 
Let $G$ be a connected $k$-regular graph with $n$ vertices 
and $m$ unoriented edges, $k>2$ and $g(G)>4$.  
Then $( {\bf U}^3)^+ $ has $2n$ eigenvalues $\lambda {}_3 $ of the form :
\[
\lambda {}_3 = \frac{1}{2} \{ \lambda {}_A ( \lambda {}^2_A -3k+4) 
\pm \sqrt{ \lambda {}^6_A -2(3k-2) \lambda {}^4_A +(13 k^2 -24k+16) \lambda {}_A -8(k-1)(k^2 -2k+2)} \} ,  
\]
where $ \lambda {}_A$ is an eigenvalues of ${\bf A} (G)$. 
The remaining $2(m-n)$ eigenvalues of $( {\bf U}^3)^+ $ are $\pm 2$ with multiplicities $m-n$. 
\end{cor}

{\bf Proof}.  By Corollary 5.3, we have 
\[
\det ( \lambda {\bf I}_{2m} - ( {\bf U}^3)^+ )=( \lambda {}^2 -4 )^{m-n} 
\prod_{ \lambda {}_{A} \in Spec( {\bf A} (G))} ( \lambda {}^2 - \lambda ( \lambda {}^3_{A} -(3k-4) \lambda {}_{A} )
+ \lambda {}^4_{A} - k^2 \lambda {}^2_{A} +2(k-1)(k^2 -2k+2)) . 
\]
Thus, solving the following equation 
\[
 \lambda {}^2 - \lambda ( \lambda {}^3_{A} -(3k-4) \lambda {}_{A} )
+ \lambda {}^4_{A} - k^2 \lambda {}^2_{A} +2(k-1)(k^2 -2k+2)=0 
\]
the result follows. 
Q.E.D.

\section{The poles and the radius of convergence of 
the 3-modified zeta function of a graph}

At first, we present the poles the 3-modified zeta function of a regular graph $G$ by means of spectra of 
the adjacency matrix ${\bf A} (G)$ of $G$.

\begin{thm}
Let $G$ be a connected $k$-regular graph with $n$ vertices 
and $m$ unoriented edges, $k>2$ and $g(G)>4$.  
Then the 3-modified zeta function ${\bf Z}_3 (G,u)$ of $G$ has $2n$ poles $\rho {}_3 $ of the form :
\[
\rho {}_3 = \frac{ \lambda {}_A ( \lambda {}^2_A -3k+4) \pm \sqrt{ \lambda {}^6_A -2(3k-2) \lambda {}^4_A 
+(13 k^2 -24k+16) \lambda {}_A -8(k-1)(k^2 -2k+2)}}{2( \lambda {}^4_{A} - k^2 \lambda {}^2_{A} +2(k-1)(k^2 -2k+2))} ,  
\]
where $ \lambda {}_A$ is an eigenvalues of ${\bf A} (G)$. 
The remaining $2(m-n)$ poles of ${\bf Z}_3 (G,u)$ are $\pm 1/2$ with multiplicities $m-n$. 
\end{thm}

{\bf Proof}.  At first, we have 
\[
{\bf Z}_3 (G,u)= \det ( {\bf I}_{2m} -u ({\bf U}^3 )^+ )
= \prod_{ \lambda {}_{3} \in Spec(({\bf U}^3 )^+ )} (1-u \lambda {}_{3} )^{-1} . 
\]
Thus, the poles of the 3-modified zeta function $ {\bf Z}_3 (G,u)$ of $G$ is given by 
\[
1/ \lambda {}_{3} , \ \lambda {}_{3} \in Spec(({\bf U}^3 )^+ ) . 
\] 

By Corollary 5.4, the $( {\bf U}^3)^+ $ has $2n$ eigenvalues $\lambda {}_3 $ of the form :
\[
\lambda {}_3 = \frac{1}{2} \{ \lambda {}_A ( \lambda {}^2_A -3k+4) 
\pm \sqrt{ \lambda {}^6_A -2(3k-2) \lambda {}^4_A +(13 k^2 -24k+16) \lambda {}_A -8(k-1)(k^2 -2k+2)} \} ,  
\]
where $ \lambda {}_A$ is an eigenvalues of ${\bf A} (G)$. 
The remaining $2(m-n)$ eigenvalues of $( {\bf U}^3)^+ $ are $\pm 2$ with multiplicities $m-n$. 
Thus, $ \pm \frac{1}{2}$ are poles of $ \tilde{{\bf Z}} (G,u)$ with multiplicity $m-n$, respectively. 
Furthermore, we have 
\[
( \frac{1}{2} \{ \lambda {}_A ( \lambda {}^2_A -3k+4) 
\pm \sqrt{ \lambda {}^6_A -2(3k-2) \lambda {}^4_A +(13 k^2 -24k+16) \lambda {}_A -8(k-1)(k^2-2k+2)} \} )^{-1}  @
\] 
\[
= \frac{ \lambda {}_A ( \lambda {}^2_A -3k+4) \pm \sqrt{ \lambda {}^6_A -2(3k-2) \lambda {}^4_A 
+(13 k^2 -24k+16) \lambda {}_A -8(k-1)(k^2 -2k+2)}}{2( \lambda {}^4_{A} - k^2 \lambda {}^2_{A} +2(k-1)(k^2 -2k+2))} . 
\]
Q.E.D.

Next, we treat the radius of convergence of the 3-modified zeta function of a graph. 
We state Perron-Frobenius Theorem on irreducible nonnegative matrices(see \cite{Ga},\cite{GR}). 

Let ${\bf A} =( a_{ij} )$ be a real square matrix. 
Then ${\bf A}$ is called {\em nonnegative} if all entries of ${\bf A}$ are nonnegative. 
Furthermore, the {\em underlying digraph} $D_{{\bf A}}$ of ${\bf A}$ has vertex set $ \{1,2, \ldots ,n \}$, 
with an arc vertex $i$ to vertex $j$ if and only if $a_{ij} \neq 0$. 
Then a nonnegative square matrix ${\bf A}$ is called {\em irreducible} if its underlying digraph is 
strongly connected. 
A digraph $D$ is {\em strongly connected} if there exists a $(x,y)$-path in $D$ for any vertices $x,y \in V(D)$.

\begin{thm}[Perron-Frobenius Theorem]
Let ${\bf A}$ be an irreducible nonnegative $n \times n$ matrix 
whose underlying digraph is strongly connected.  
Furthermore, let ${\bf u} = {\bf A} {\bf 1} $, where 
${\bf 1} $ is the vector with all one. 
Set ${\bf u} = {}^t (u_1 , \ldots , u_n )$.  
Then 
\begin{enumerate}
\item  ${\bf A}$ has at least one positive eigenvalue. 
The maximal positive eigenvalue $\alpha $ is simple. 

\item $| \lambda | \leq \alpha $ for any eigenvalue $ \lambda $ of ${\bf A}$. 

\item 
\[
{\rm min}_{1 \leq i \leq n} u_i \leq \alpha \leq {\rm max}_{1 \leq i \leq n} u_i . 
\]
\end{enumerate} 
\end{thm}

In Theorem 6.2, $\alpha $ is called the {\em spectral radius} of ${\bf A}$. 

Now, we consider the strongly connectivity of the underlying digraph $D_{({\bf U}^3 ) ^+} $ 
of $({\bf U}^3 ) ^+$. 
Let $G$ be a connected graph with $n$ vertices and $m$ edges. 
For a path $P=(e_1, \ldots , e_r)$, we say that $P$ is called an {\em $(e_1 ,e_r)$-path}.

\begin{pro}
Let $G$ be a connected $k$-regular graph with $n$ vertices and $m$ edges, and $X= D_{({\bf U}^3 ) ^+}$ 
the underlying digraph of $({\bf U}^3 ) ^+ $.  
Suppose that $k>2$ and $g(G)>4$. 
Then $X$ is strongly connected. 
\end{pro}

{\bf Proof}.  Let $e,f \in D(G)$. 
By the definition of 3-paths, there exists an $(e,f)$-3-path in $G$ if and only if 
there exists an $(e,f)$-path in $X$.

Now, let $e,f \in D(G)$ and $P=(e_1 , e_2 , \ldots , e_r )$ an $(e,f)$-path in $G$, where $e=e_1 $ and $f=e_r $. 
Suppose that $P$ is reduced.  
If $|P|=4k$, then $P$ is an $(e,f)$-3-path without 3-backtracking, and so 
there exists an $(e,f)$-path in $X$. 

Next, let $P|=4k+j(j=1,2,3)$. 
Since $k>2$, there exists a cycle $C=(e_1 , e_2 , \ldots , e_r , e_{r+1} , \ldots , e_q )$ 
in $G$ such that $P \subset C$. 
Without of generality, we may assume that $C$ has no backtracking. 
Let 
\[
Q=(e_1 , e_2 , \ldots , e_r , \ldots , e_{r+1} , \ldots , e_{r+4-j} ) \subset C . 
\]
Then $|Q|=4k+4$. 

Let 
\[
R=(e_1 , e_4 , e_8 , \ldots , e_{r+4-j}, e_{r+3-j} , \ldots , e_{r} ) , 
\]
where $(e_{r+4-j}, e_{r+3-j} , \ldots , e_{r} )$ is a sequence of 3-backtrackings. 
Then $R$ is an $(e,f)$-3-path in $G$. 
Thus, there exists an $(e,f)$-path in $X$. 
Therefore, $X$ is strongly connected. 
Q.E.D.

Thus, we obtain the following result.

\begin{cor}

Let $G$ be a connected $k$-regular graph with $n$ vertices and $m$ edges, $k>2$ and $g(G)>4$. 
Furthermore, let $\rho $ be the radius of convergence of the 3-modified zeta function $ {\bf Z}_3 (G,u)$ of $G$. 
Then the following holds for  
\[ 
\rho = \frac{1}{k^3 +k-1} . 
\]
\end{cor}

{\bf Proof}.  By Proposition 4.1, we have 
\[
({\bf U}^3)^+ =( {\bf U}^+ )^3 + {}^T {\bf U}^+ . 
\] 
Let $e \in D(G)$. 
Then we consider the sum $r_e$ of the $e$-row of $( {\bf U}^3 )^+$. 
By the definition of 3-arcs and 3-backtrackings, $r_e$ is equal to 
the number of 3-arcs $[e,f]$ and 3-backtrackings $[f,e]$. 
Thus, 
\[
r_e = k^3 +k-1 .  
\]
By Theorem 6.2, it follows that 
\[
\rho = \frac{1}{k^3 +k-1} . 
\]
Q.E.D.

\
\par\noindent
\section*{Acknowledgments}
YuH's work was supported in part 
by JSPS Grant-in-Aid for Scientific Research (C)~20540113, 25400208 
and (B)~24340031. 
NK and IS also acknowledge financial supports of 
the Grant-in-Aid for Scientific Research (C) 
from Japan Society for the Promotion of Science 
(Grant No.~24540116 and No.~23540176, respectively).
ES thanks to the financial support of
the Grant-in-Aid for Young Scientists (B) of Japan Society for the
Promotion of Science (Grant No. 25800088). 
\par

\begin{small}
\bibliographystyle{jplain}

\end{small}

\end{document}